\documentclass{ifacconf}

\usepackage{graphicx} 
\usepackage{epstopdf} 
\usepackage{psfrag} 
\usepackage{amsmath} 
\usepackage{amssymb}  
\usepackage{algorithm}
\usepackage{algorithmic}
\usepackage{thmtools,dsfont}
\usepackage{color}
\usepackage{tabu}

\usepackage{enumitem} \setlist[itemize,1]{leftmargin=\dimexpr 24pt}
\usepackage{cite}
\usepackage{microtype,enumerate,url}
\usepackage{BOONDOX-ds}
\usepackage{tikz} 
\usetikzlibrary{positioning} 
\usepackage[fulladjust]{marginnote}
\usepackage{pgfplots} %
\usepgfplotslibrary{fillbetween}
\pgfplotsset{compat=1.17} 
\usepackage{pgfplotstable} 
\usepackage{xstring} 
\usepackage{booktabs} 
\usepackage{colortbl}
\usepackage[skip=0pt]{caption}
\usepackage{mathdots} 
\usepackage{tcolorbox}
\usepackage[normalem]{ulem} 
\allowdisplaybreaks 
\usepackage{tikz-cd}
\usepackage[switch, modulo]{lineno}

\graphicspath{{./img/}}

\renewcommand{\theenumi}{(\roman{enumi})}

\newtheorem{definition}{\bf Definition}[section]

\newcommand{\beq}{\begin{equation}}
\newcommand{\eeq}{\end{equation}}
\newcommand {\bseq}{\begin{subequations}}
\newcommand {\eseq}{\end{subequations}}

\newcommand {\R}{\mathbb{R}} 	

\newcommand{\Image}{\operatorname{im}} 
\newcommand{\transpose}{\mathsf{T}} 

\newcommand{\norm}[1]{\left\lVert#1\right\rVert}

\newcommand{\trace}{\operatorname{trace}}

\newcommand\Grass[2]{\operatorname{Gr}(#1,  #2) }
\newcommand\St[2]{\operatorname{St}(#1,  #2) }

\definecolor{green_color_blind}{RGB}{102,194,165}
\definecolor{orange_color_blind}{RGB}{252,141,98}
\definecolor{blue_color_blind}{RGB}{141,160,203} 
\definecolor{jeremy_orange}{RGB}{230, 159, 0}
\definecolor{jeremy_lightblue}{RGB}{86, 180, 233}
\definecolor{jeremy_green}{RGB}{0, 158, 115}
\definecolor{jeremy_yellow}{RGB}{240, 228, 66}
\definecolor{jeremy_blue}{RGB}{0, 114, 178}
\definecolor{jeremy_vermillion}{RGB}{213, 94, 0}
\definecolor{jeremy_pink}{RGB}{204, 121, 167}

\usepackage{natbib}        

\begin{document}
\begin{frontmatter}

\title{Geometrically robust least squares through manifold optimization}

\author[First]{Jeremy Coulson} 
\author[Second]{Alberto Padoan} 
\author[Third]{Cyrus Mostajeran} 

\address[First]{University of Wisconsin–Madison (e-mail: jeremy.coulson@wisc.edu)}
\address[Second]{ETH Z\"urich (e-mail: apadoan@ethz.ch).}
\address[Third]{School of Physical and Mathematical Sciences, Nanyang Technological University, Singapore (email: cyrussam.mostajeran@ntu.edu.sg).}

\begin{abstract}
This paper presents a methodology for solving a geometrically robust least squares problem, which arises in various applications where the model is subject to geometric constraints. The problem is formulated as a minimax optimization problem on a product manifold, where one variable is constrained to a ball describing uncertainty. To handle the constraint, an exact penalty method is applied. A first-order gradient descent ascent algorithm is proposed to solve the problem, and its convergence properties are illustrated by an example. The proposed method offers a robust approach to solving a wide range of problems arising in signal processing and data-driven control.
\end{abstract}
\end{frontmatter}

\section{Introduction}\label{sec:introduction}

The \textit{least squares problem} is a cornerstone of science and engineering, classically defined by the optimization problem
\begin{equation} \label{eq:least_squares}
\underset{x \in \mathbb{R}^{k}}{\textup{min}}\;
\|A x - b\|_2^2,
\end{equation}
where $x \in \R^k$ is an unknown coefficient vector, while $A \in \R^{n \times k}$ and $b \in \R^n$ are given data which usually represent a \textit{linear model} (often termed the \textit{regressor matrix} or the \textit{design matrix}) and a vector of observations, respectively. 
The historical origins of the least squares problem can be traced back to its early applications in celestial mechanics~\citep{gauss1809,gauss1857} and it has since become indispensable across various scientific and engineering fields. In control theory, least squares problems are ubiquitous and play a critical role in numerous applications, including  identification, filtering, and optimal control, among others.

The least squares problem~\eqref{eq:least_squares} is an essential tool to approximate solutions of the linear system ${Ax = b}$. By minimizing the residual $\Delta b = Ax-b$, it determines a solution that closely matches the observations in the Euclidean $2$-norm. 
The \textit{total least squares problem} discussed in~\citep{golub1980analysis} extends this paradigm 
to address the presence of \textit{perturbations} (errors or uncertainties) in both the dependent and independent variables.  This extension involves minimizing residuals $(\Delta A, \Delta b)$ in the observations $b$ and in the matrix $A$, thereby refining the accuracy of the model. Given that accuracy is the primary aim of both least squares and total least squares, it is not surprising that both solutions may exhibit sensitive behavior to perturbations in the data matrices $(A,b)$~\citep{golub2013matrix}.

Sensitivity issues and robustness to perturbations in the data $(A,b)$ can be effectively mitigated in different ways.
One strategy is to consider the \textit{robust} least squares problem
\begin{equation} \label{eq:least_squares_robust_Tikhonov}
\underset{x \in \R^{k}}{\textup{min}}\;
\underset{A \in \mathbb{B}^{\scriptscriptstyle \text{F}}_{\rho}(\hat{A})}{\textup{max}}\;
\|A x - b\|_2^2 ,
\end{equation}
where
$\mathbb{B}^{\scriptscriptstyle \text{F}}_{\rho}(\hat{A}) = \{A \in \mathbb{R}^{n\times k} \, |\, \|A- \hat{A}\|_{\text{F}} \le \rho \}$
is the ball centered at $\hat{A}$ with radius $\rho$ defined by the Frobenius norm.
This approach is intimately connected to (Tikhonov) \textit{ regularization}~\citep{golub2013matrix}.

The perturbation model defined implicitly by~\eqref{eq:least_squares_robust_Tikhonov} allows the perturbed data matrix $A$ to vary in a neighbourhood of the matrix $\hat{A}$. The parameter $\rho$ controls the size of the neighbourhood and can be interpreted as a penalty that promotes solutions whose norms are small, thereby trading accuracy for robustness at the expense of introducing bias.
Consequently, the formulation \eqref{eq:least_squares_robust_Tikhonov} serves a double purpose: it enhances the robustness of the original problem and it incorporates prior knowledge into the model.

Building on these ideas, alternative perturbation models have been explored, yielding different robust versions of the least squares problem. For example, ~\citep{el1997robust} considered the robust least squares problem 
\begin{equation} \label{eq:least_squares_regularized_El_Ghaoui}
\underset{x \in \R^{k}}{\textup{min}}\;
\underset{[\, A \ b\,]  \in 
\mathbb{B}^{\scriptscriptstyle \text{F}}_{\rho}
([\, \hat{A} \ \hat{b}\,])}{\textup{max}}\; \|A x - b\|_2^2 .
\end{equation}
This approach is motivated by scenarios where the exact data $(A, b)$ are unknown, but belong to a family of matrices $(\hat{A}+\Delta A, \hat{b}+\Delta b)$ and the residual $[\,\Delta A \ \Delta b \,]$ lies in a norm-bounded matrix ball. The paper also considers a structured variant of~\eqref{eq:least_squares_regularized_El_Ghaoui} by considering matrix box constraints of given radii \citep{el1997robust}.

In this extended abstract, we introduce a robust least squares problem that accounts for the geometric nature of perturbations found in diverse instances of the problem. Specifically, we consider the optimization problem
\begin{equation} \label{eq:least_squares_robust_geometric}
\underset{x \in \R^{n}}{\textup{min}}\;
\underset{S  \in \mathbb{B}^d_{\rho}(\hat{S})}{\textup{max}}\; \|P_S x - b\|_2^2 ,
\end{equation}
where $P_{S}$ is the orthogonal projection matrix onto the $k$-dimensional subspace $S$ of $\R^{n}$ and $\mathbb{B}_{\rho}^d(\hat{S})$ is the ball centered at $\hat{S}$ with radius $\rho$ defined by the metric $d$ on the \textit{Grassmannian} $\Grass{k}{n}$, i.e., the set of all $k$-dimensional subspaces in $\R^n$ endowed with the structure of a smooth Riemannian manifold. This approach is motivated by a diverse range of applications where the linear model $A$ is a matrix representation of a subspace subject to bounded perturbations (due to uncertainty or approximations errors), quantified naturally in terms of distances between subspaces. In what follows, we present some examples to motivate~\eqref{eq:least_squares_robust_geometric}.

\textit{Subspace tracking:} In signal processing, communications, and computer vision, the problem of tracking a time-varying $k$-dimensional subspace $S_t$ in $\mathbb{R}^{n}$ is common and arises in various contexts~\citep{delmas2010subspace}. Subspace tracking algorithms aim to provide an estimate $\hat{S}_t$ of the true subspace from a sequence of possibly noisy samples in $S_t$. The projection matrix $P_{\hat{S}_t}$ onto the subspace estimate $\hat{S}_t$ is then used as a linear model, e.g., for foreground-background separation in noisy free-motion camera videos~\citep{he2012incremental} or for Direction of Arrival (DOA) tracking in blind channel estimation and equalization~\citep{delmas2010subspace}. Subspace tracking algorithms often guarantee local~\citep{onlineidentification2010} or global~\citep{globalconvergence2016} convergence with respect to a given metric on the Grassmannian (e.g., the chordal metric in~\citep{globalconvergence2016}). 
Consequently, quantifying the distance between the estimate and the true subspace after a finite number of iterations of a subspace tracking algorithm is naturally expressed by a metric ball. This, in turn, leads to the robust least squares problem~\eqref{eq:least_squares_robust_geometric} when using the projection matrix $P_{\hat{S}_t}$ as a linear model, which needs to be robustified up to the maximum distance within the metric ball.

\textit{Data-driven control:}
The availability of large datasets coupled with unprecedented storage and computational power has recently reignited interest in direct data-driven control methods, which aim to infer optimal decisions directly from measured data (bypassing system identification).  At the heart of this emerging trend lies the behavioral approach to system theory~\citep{willems2007behavioral} and a seminal result by Willems and his collaborators~\citep{willems2005note}, commonly referred to as the \textit{fundamental lemma}. The lemma establishes that finite-horizon behaviors of Linear Time-Invariant (LTI) systems can be represented as  images of raw data matrices. As a result, various data-driven modeling, estimation, filtering, and control problems can be formulated as weighted or constrained least squares problems~\citep{markovsky2021behavioral}. Moreover, being finite-dimensional subspaces,  finite-horizon LTI behaviors can be identified with points on the Grassmannian $\Grass{k}{n}$ and uncertainty can be naturally quantified using Grassmannian (subspace) metrics. This approach has demonstrated its effectiveness in data-driven mode recognition and control applications and shows promise to open new avenues in adaptive control~\citep{padoan2019behavioral}.

The key contribution of this paper is the formulation of a methodology for solving the geometrically robust least squares problem \eqref{eq:least_squares_robust_geometric}, grounded in the framework of optimization on manifolds \citep{absil2008optimization, boumal2023introduction} and enriched by recent progress in minimax optimization techniques~\citep{han2023nonconvexnonconcave}.
Due to the presence of a (ball) constraint in the second variable, problem \eqref{eq:least_squares_robust_geometric}  cannot be solved directly using the framework presented in \citep{han2023nonconvexnonconcave}. We resolve this issue using an exact penalty method. In particular, we augment the original cost function with a penalty term to enforce the ball constraint. We then apply a smoothing technique~\citep{liu2020simple} to transform the resulting non-smooth penalized cost function into a smooth one. Finally, we introduce a first-order gradient descent ascent algorithm, motivated and inspired by the TSRGDA algorithm detailed in \citep{han2023nonconvexnonconcave}. The convergence properties of the algorithm are demonstrated by means of an illustrative example. 

\section{Geometrically robust least squares}\label{sec:problem}

We adopt optimization on manifolds as our framework \citep{absil2008optimization, boumal2023introduction}, given that one optimization variable of the geometrically robust least squares problem~\eqref{eq:least_squares_robust_geometric} naturally lies in a nonlinear space. Our goal is to turn~\eqref{eq:least_squares_robust_geometric} into a minimax optimization problem 
\beq \label{eq:minimax-general}
\displaystyle \min_{x \in \mathcal{X}} \max_{y \in \mathcal{Y}} f(x,y) ,
\eeq
where $\mathcal{X}$ and $\mathcal{Y}$ are subsets of smooth finite-dimensional Riemannian manifolds $\mathcal{M}$ and $\mathcal{N}$, respectively. 

Problem~\eqref{eq:least_squares_robust_geometric} fits into this minimax framework with the cost function given by
\begin{equation} \label{eq:cost}
f(x,y)  = \|P_{y} x - b\|_2^2  .
\end{equation}
The manifold $\mathcal{M}$ is the Euclidean space $\R^n$ endowed with the Riemannian metric defined by the standard inner product. 
The manifold $\mathcal{N}$ is the Grassmannian $\Grass{k}{n}$ endowed with the Riemannian metric
defined in~\citep[Section 9.16]{boumal2023introduction}.
Furthermore, $\mathcal{X}=\R^n$ and 
$\mathcal{Y}=\mathbb{B}_{\rho}^d(\hat{y})$ is the ball centered at $\hat{y}$ with radius $\rho$ defined by the metric $d$ (not necessarily the Riemannian distance).

Since the cost function~\eqref{eq:cost} is smooth and convex in $x$, but non-concave in $y$, we assume that the cost function $f: \mathcal{M} \times \mathcal{N} \to \R$ in~\eqref{eq:minimax-general} is smooth, but in general non-convex and non-concave (in both Euclidean and geodesic senses). As a result, the existence of saddle points is not guaranteed~\citep{jin2020local}. Consequently, we adopt the notion of solution defined by a local minimax point from~\citep{han2023nonconvexnonconcave}, where $d_{\mathcal{M}}$ denotes the Riemannian distance in $\mathcal{M}$ (i.e., the length of shortest path connecting $z_1,z_2 \in  \mathcal{M}$).
\begin{definition}[Local minimax point]\label{def:minimax}
    A point $(x_*,y_*)$ is called a \emph{local minimax point} of $f:\mathcal{X}\times \mathcal{Y} \to \R$ if there exists a constant $\delta_0>0$ and a function $h$ satisfying $\lim_{\delta \to 0} h(\delta)=0$ for any $\delta\in(0,\delta_0]$ and any $(x,y)\in \mathcal{X}\times \mathcal{Y}$ satisfying $d_{\mathcal{M}}(x,x_*) \leq \delta$ and $d_{\mathcal{N}}(y,y_*)\leq\delta$, we have $f(x_*,y)\leq f(x_*, y_*) \leq \textup{max}_{y^{\prime}:d_{\mathcal{N}}(y',y_*)\leq h(\delta)} f(x, y')$.
\end{definition}
A necessary condition for $(x_*,y_*)$ to be a local minimax point is that the Riemannian (partial) gradients denoted by $\textup{grad}_x f(x_*, y_*)$ and $\textup{grad}_y f(x_*, y_*)$ are zero~\citep[Proposition 3]{han2023nonconvexnonconcave}.
Sufficient conditions in terms of the Hessian can be found in~\citep{han2023nonconvexnonconcave}.

Existing algorithms for solving unconstrained minimax optimization problems on manifolds \citep{han2023nonconvexnonconcave} do not explicitly accommodate (metric ball) constraints. To overcome this limitation, we use an exact penalty method \citep{liu2020simple}. Introducing a penalty weight ${\lambda \geq 0}$, we relax the original problem and consider the unconstrained minimax optimization problem 
\beq \label{eq:minimax-relaxed}
\displaystyle \min_{x \in \mathcal{M}} \max_{y \in \mathcal{N}} 
f(x,y) +\lambda\min\{0,\rho^2-d(y,\hat{y})^2\}.
\eeq 
Given that this lifting results in a nonsmooth objective, we employ a smoothing technique inspired by \citep{liu2020simple} to obtain a smooth relaxation of the minimax optimization problem~\eqref{eq:minimax-relaxed}, obtaining
\beq \label{eq:minimax-relaxed-smoothed}
\displaystyle \min_{x \in \mathcal{M}} \max_{y \in \mathcal{N}} 
f(x,y) -\lambda u 
\log\left(1+\textup{exp}\left(\tfrac{d(y,\hat{y})^2 - \rho^2}{u}\right)\right)
\eeq 
where ${u>0}$ serves as the smoothing parameter.
For the specific case of the geometrically robust least squares problem~\eqref{eq:least_squares_robust_geometric}, these transformations yield
\begin{equation} \label{eq:least_squares_robust_smooth}
    \underset{x \in \R^{n}}{\textup{min}}\;
\underset{ y \in \Grass{k}{n}}{\textup{max}}\; \|P_{y} x - b\|_2^2 -\lambda u \log\left(1+\textup{exp}\left(\tfrac{d(y,\hat{y})^2 - \rho^2}{u}\right)\right).
\end{equation}

The minimax optimization problem can then be solved using any first- or second-order algorithm presented in~\citep{han2023nonconvexnonconcave}. For simplicity, we opt for the \textit{timescale-separated Riemannian gradient descent ascent} (TSRGDA) algorithm, which is based on the update rule 
\begin{equation} \label{eq:update}
(x_{k+1},y_{k+1}) = \textup{Exp}_{(x_k,y_k)}\begin{pmatrix}
   -\eta_x \textup{grad}_x f(x_k,y_k) \\
   \eta_y \textup{grad}_y f(x_k,y_k) 
   \end{pmatrix} ,
\end{equation}
where $\textup{Exp}_{z}$ denotes the exponential map of the product manifold $\mathcal{M}\times\mathcal{N}$, and $\eta_x,\eta_y>0$ are step sizes.

In our numerical implementations, each point $y \in \Grass{k}{n}$ is identified with an equivalence class represented by the image of an orthonormal matrix $Y\in \R^{n\times k}$,  treated as an element of the \emph{Stiefel manifold} 
 \[
 \textup{St}(k,n):=\{Y\in\R^{n\times k}\,|\,Y^\top Y=I\}.
 \]
The Stiefel manifold is viewed as a Riemannian submanifold of the Euclidean space $\R^{n\times k}$, endowed with Riemannian metric inherited by the standard Euclidean inner product
\[
\langle Z_1,\, Z_2\rangle = \trace(Z_1^{\top}Z_2).
\]
This approach is standard and the interested reader is referred to~\citep{absil2008optimization} for further detail.

The metric $d$ used to describe the ball constraint in~\eqref{eq:least_squares_robust_smooth} is chosen as the \textit{chordal distance} between subspaces
\[
d_{\Grass{k}{n}}^{\kappa}(y_1, y_2)=  \left(\sum_{i=1}^k (\sin\theta_i)^2\right)^{1/2} ,   
\]
where $\theta_{i}$ denotes the $i$-th principal angle between the subspaces $y_1\in \Grass{k}{n}$ and $y_2\in \Grass{k}{n}$ (see~\citep{golub2013matrix} for definition of principal angles). Given $Y_1,Y_2\in \St{k}{n}$ such that  $y_1=\Image(Y_1)$ and $y_2=\Image(Y_2)$, the chordal distance can be computed as
\beq
d_{\Grass{k}{n}}^{\kappa}(y_1, y_2) = \tfrac{1}{\sqrt 2} \norm{ Y_1Y_1^{\transpose} - Y_2Y_2^{\transpose}}_F  .
\eeq
This distance does not correspond to the standard Riemannian distance, but it  satisfies the four metric requirements while being more numerically stable than the geodesic distance. Thus, it is preferred over the Riemannian distance here (and in most applications).

\section{Example}
Consider problem~\eqref{eq:least_squares_robust_smooth} with $n=2$, $k=1$,
\[
\hat{y}=\Image\begin{pmatrix} 1 \\ 0\end{pmatrix},\quad b=\begin{pmatrix}\cos(\tfrac{\pi}{16}) \\ \sin(\tfrac{\pi}{16})\end{pmatrix},
\]
subspace metric $d=d_{\Grass{k}{n}}^{\kappa}$, radius $\rho=\sin(\tfrac{\pi}{8})$, smoothing parameter $u=0.01$, and penalty $\lambda=70$. Starting from a random initial point $(x_0,y_0)\in\R^2\times \Grass{1}{2}$, we apply update rule~\eqref{eq:update} with step sizes $\eta_x=0.01$, and $\eta_y=0.1$. This results in iterates $x_k$ and $y_k$ that converge to $x_*$ and $y_*$. The projected iterates $P_{y_k}x_k$, $P_{y_*}x_*$, and $y_*$ are illustrated in Fig.~\ref{fig:geometry}. We see from Fig.~\ref{fig:gradients} that the norm of the gradient of the objective function in~\eqref{eq:least_squares_robust_smooth} converges to zero. Thus, $(x_*,y_*)$ satisfies the necessary condition for being a minimax point (see the discussion after Definition~\ref{def:minimax}).
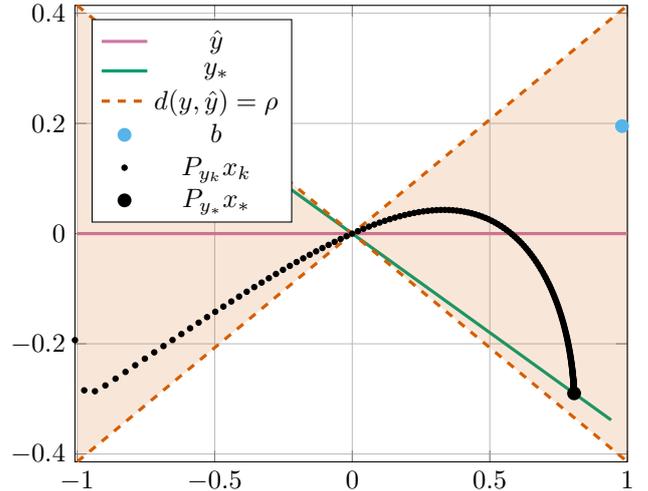
\begin{figure}[h]
\centering
\begin{tikzpicture}[black, every tick label/.append style={font=\normalsize}]
\begin{axis}[
enlargelimits=false,
width = \columnwidth,
grid,
legend pos=north west,
legend columns=1, 
]
\addplot [name path=S, jeremy_pink, very thick]  table [x index = {0}, y index = {1},  col sep=comma]{S.csv};
\addlegendentry{$\hat{y}$}
\addplot [jeremy_green, very thick]  table [x index = {0}, y index = {1},  col sep=comma]{S_star.csv};
\addlegendentry{$y_*$}
\addplot [name path=A, jeremy_vermillion, very thick, dashed]  table [x index = {0}, y index = {1},  col sep=comma]{ball_boundary_upper.csv};
\addlegendentry{$d(y,\hat{y})=\rho$}
\addplot [name path=B, jeremy_vermillion, very thick, dashed, forget plot]  table [x index = {0}, y index = {1},  col sep=comma]{ball_boundary_lower.csv};
\tikzfillbetween[of=B and S]{jeremy_vermillion, opacity=0.15};
\tikzfillbetween[of=S and A]{jeremy_vermillion, opacity=0.15};
\addplot [jeremy_lightblue, very thick, mark=*, only marks]  table [x index = {0}, y index = {1},  col sep=comma]{b.csv};
\addlegendentry{$b$}
\addplot [black, mark=*, mark size=1pt, only marks]  table [x index = {0}, y index = {1},  col sep=comma]{x_iterates.csv};
\addlegendentry{$P_{y_k}x_k$}
\addplot [black, mark=*, very thick, only marks]  table [x index = {0}, y index = {1},  col sep=comma]{x_star.csv};
\addlegendentry{$P_{y_*}x_*$}

\end{axis}
\end{tikzpicture}
\caption{
The points $y_*$ and $P_{y_*}x_*$ obtained from iteratively applying update rule~\eqref{eq:update} are plotted as a green solid line and large black dot, respectively. The projected iterates $P_{y_k}x_k$ are plotted as small black dots. The boundaries of $\mathbb{B}^d_{\rho}(\hat{y})$ are depicted with dashed lines and its interior is shaded.
}
\label{fig:geometry}
\end{figure}%
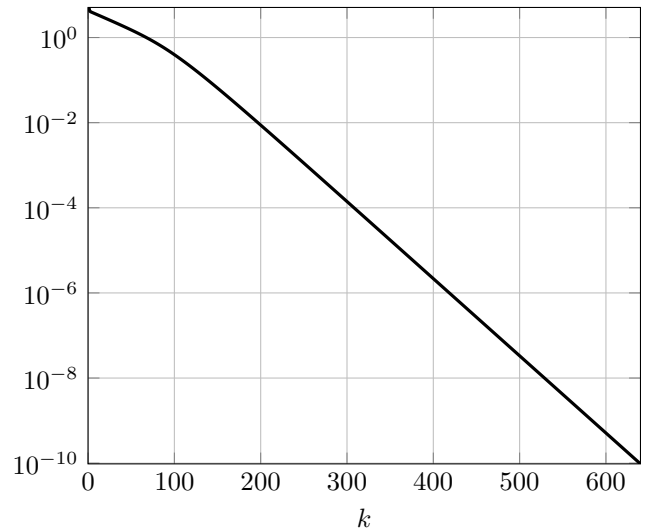
\begin{figure}[h]
\centering
\begin{tikzpicture}[black, every tick label/.append style={font=\normalsize}]
\begin{axis}[
enlargelimits=false,
width = \columnwidth,
ymode=log,
xlabel near ticks,
xlabel={$k$},
ylabel near ticks,
grid,
legend pos=north east
]
\addplot [black, very thick]  table [x index = {0}, y index = {1},  col sep=comma]{gradient_norm.csv};
\end{axis}
\end{tikzpicture}
\caption{
Frobenius norm of the gradient of the objective function in~\eqref{eq:least_squares_robust_smooth} at the iterates of update rule~\eqref{eq:update}.
}
\label{fig:gradients}
\end{figure}%

\section{Conclusion}\label{sec:conclusion}
We have introduced a robust least squares problem tailored to accommodate geometric perturbations commonly encountered in various scenarios. Leveraging the framework of optimization on manifolds and recent advances in minimax optimization, we have addressed the problem through an exact penalty method and a smoothing technique. We proposed a first-order gradient descent ascent algorithm to solve the problem, illustrating its convergence properties through a numerical example. Our methodology not only provides a new robust least squares formulation but also expands its applicability to a diverse range of problems where the geometric nature of the linear model is essential. 


\begin{thebibliography}{18}
\providecommand{\natexlab}[1]{#1}
\providecommand{\url}[1]{\texttt{#1}}
\providecommand{\urlprefix}{URL }
\expandafter\ifx\csname urlstyle\endcsname\relax
  \providecommand{\doi}[1]{doi:\discretionary{}{}{}#1}\else
  \providecommand{\doi}{doi:\discretionary{}{}{}\begingroup
  \urlstyle{rm}\Url}\fi

\bibitem[{Absil et~al.(2008)Absil, Mahony, and
  Sepulchre}]{absil2008optimization}
Absil, P.A., Mahony, R., and Sepulchre, R. (2008).
\newblock \emph{Optimization algorithms on matrix manifolds}.
\newblock Princeton University Press, Princeton.

\bibitem[{Balzano et~al.(2010)Balzano, Nowak, and
  Recht}]{onlineidentification2010}
Balzano, L., Nowak, R., and Recht, B. (2010).
\newblock Online identification and tracking of subspaces from highly
  incomplete information.
\newblock In \emph{2010 48th Annual Allerton Conference on Communication,
  Control, and Computing (Allerton)}, 704--711. IEEE.

\bibitem[{Boumal(2023)}]{boumal2023introduction}
Boumal, N. (2023).
\newblock \emph{An introduction to optimization on smooth manifolds}.
\newblock Cambridge University Press, Cambridge.

\bibitem[{Delmas(2010)}]{delmas2010subspace}
Delmas, J.P. (2010).
\newblock Subspace tracking for signal processing.

\bibitem[{El~Ghaoui and Lebret(1997)}]{el1997robust}
El~Ghaoui, L. and Lebret, H. (1997).
\newblock Robust solutions to least-squares problems with uncertain data.
\newblock \emph{SIAM Journal on Matrix Analysis and Applications}, 18(4),
  1035--1064.

\bibitem[{Gauss(1809)}]{gauss1809}
Gauss, C.F. (1809).
\newblock \emph{Theoria motus corporum coelestium in sectionibus conicis solem
  ambientum}.
\newblock Hamburg.
\newblock \url{https://archive.org/details/bub_gb_ORUOAAAAQAAJ}.

\bibitem[{Gauss(1857)}]{gauss1857}
Gauss, C.F. (1857).
\newblock \emph{Theory of the motion of the heavenly bodies moving about the
  sun in conic sections: a translation of Gauss’s “Theoria motus.”}.
\newblock \url{https://archive.org/details/motionofheavenly00gausrich}.

\bibitem[{Golub and Van~Loan(2013)}]{golub2013matrix}
Golub, G.H. and Van~Loan, C.F. (2013).
\newblock \emph{Matrix computations (4th Ed.)}.
\newblock Johns Hopkins Univ. Press, Baltimore, MD, USA.

\bibitem[{Golub and Van~Loan(1980)}]{golub1980analysis}
Golub, G.H. and Van~Loan, C.F. (1980).
\newblock An analysis of the total least squares problem.
\newblock \emph{SIAM journal on numerical analysis}, 17(6), 883--893.

\bibitem[{Han et~al.(2023)Han, Mishra, Jawanpuria, and
  Gao}]{han2023nonconvexnonconcave}
Han, A., Mishra, B., Jawanpuria, P., and Gao, J. (2023).
\newblock Nonconvex-nonconcave min-max optimization on {R}iemannian manifolds.
\newblock \emph{Transactions on Machine Learning Research}.

\bibitem[{He et~al.(2012)He, Balzano, and Szlam}]{he2012incremental}
He, J., Balzano, L., and Szlam, A. (2012).
\newblock {Incremental gradient on the {G}rassmannian for online foreground and
  background separation in subsampled video}.
\newblock In \emph{2012 IEEE Conference on Computer Vision and Pattern
  Recognition}, 1568--1575. IEEE.

\bibitem[{Jin et~al.(2020)Jin, Netrapalli, and Jordan}]{jin2020local}
Jin, C., Netrapalli, P., and Jordan, M. (2020).
\newblock What is local optimality in nonconvex-nonconcave minimax
  optimization?
\newblock In \emph{International conference on machine learning}, 4880--4889.
  PMLR.

\bibitem[{Liu and Boumal(2020)}]{liu2020simple}
Liu, C. and Boumal, N. (2020).
\newblock Simple algorithms for optimization on {R}iemannian manifolds with
  constraints.
\newblock \emph{Applied Mathematics \& Optimization}, 82, 949--981.

\bibitem[{Markovsky and D{\"o}rfler(2021)}]{markovsky2021behavioral}
Markovsky, I. and D{\"o}rfler, F. (2021).
\newblock Behavioral systems theory in data-driven analysis, signal processing,
  and control.
\newblock \emph{Ann. Rev. Control}, 52, 42--64.

\bibitem[{Padoan et~al.(2022)Padoan, Coulson, van Waarde, Lygeros, and
  D\"orfler}]{padoan2019behavioral}
Padoan, A., Coulson, J., van Waarde, H.J., Lygeros, J., and D\"orfler, F.
  (2022).
\newblock Behavioral uncertainty quantification for data-driven control.
\newblock In \emph{\CDC{61st}}, 4726--4731. Cancun, Mexico.

\bibitem[{Willems(2007)}]{willems2007behavioral}
Willems, J.C. (2007).
\newblock The behavioral approach to open and interconnected systems.
\newblock \emph{IEEE Control Syst. Mag.}, 27(6), 46--99.

\bibitem[{Willems et~al.(2005)Willems, Rapisarda, Markovsky, and
  De~Moor}]{willems2005note}
Willems, J.C., Rapisarda, P., Markovsky, I., and De~Moor, B.L.M. (2005).
\newblock A note on persistency of excitation.
\newblock \emph{\SCL}, 54(4), 325--329.

\bibitem[{Zhang and Balzano(2016)}]{globalconvergence2016}
Zhang, D. and Balzano, L. (2016).
\newblock Global convergence of a {G}rassmannian gradient descent algorithm for
  subspace estimation.
\newblock In \emph{Artificial Intelligence and Statistics}, 1460--1468. PMLR.

\end{thebibliography}

\newcommand{\TAC}{\textit{{IEEE} Trans. Autom.
  Control}}\newcommand{\TCST}{\textit{{IEEE} Trans. Syst.
  Tech.}}\newcommand{\TIT}{\textit{{IEEE} Trans. Inform.
  Theory}}\newcommand{\TSP}{\textit{{IEEE} Trans. Sign.
  Proc.}}\newcommand{\SCL}{\textit{Syst. Control
  Lett.}}\newcommand{\IJC}{\textit{Int. J.
  Control}}\newcommand{\EJC}{\textit{Eur. J.
  Control}}\newcommand{\ACC}{\textit{Proc. Amer. Control
  Conf.}}\newcommand{\ECC}{\textit{Proc. Eur. Control
  Conf.}}\newcommand{\CDC}[1]{\textit{Proc. {#1} Conf. Decision
  Control}}\newcommand{\CDCs}[1]{\textit{{#1} Conf. Decision
  Control}}\newcommand{\IFAC}[1]{\textit{Proc. {#1} IFAC World
  Congr.}}\newcommand{\SIAM}{\textit{SIAM J. Control
  Optim.}}\newcommand{\MCSS}{\textit{Math. Control, Sign.
  Syst.}}\newcommand{\NOLCOS}[1]{\textit{Proc. {#1} IFAC Symp. Nonlinear
  Control Syst.}}\newcommand{\MTNS}[1]{\textit{ Proc. {#1} Math. Symp. Netw.
  Syst.}}

\end{document}